\documentclass[doublecol,figures]{epl2} 

\usepackage[utf8]{inputenc}
\usepackage{amsmath}
\usepackage{amsfonts}
\usepackage{amssymb}
\usepackage{enumerate}  
\usepackage{hyperref}    
\usepackage[dvipsnames]{xcolor}
\usepackage[switch]{lineno}

\newcommand{\meantr}{\langle t_r \rangle}

\title{Stochastic gel-shatter cycles in coalescence-fragmentation models}
\shorttitle{Gel-shatter cycles} 

\author{Brennen T. Fagan\inst{1} \and Niall J. MacKay\inst{1} \and Dmitri O. Pushkin\inst{1} \and A. Jamie Wood\inst{1,2}}
\shortauthor{Fagan \etal}

\institute{                    
  \inst{1} Department of Mathematics, University of York, York, YO10 5DD, UK\\
  \inst{2} Department of Biology, University of York, York, YO10 5DD, UK
}

\pacs{36.40.Qv}{Stability and fragmentation of clusters} 
\pacs{82.40.Bj}{Oscillations, chaos, and bifurcations} 
\pacs{47.57.eb}{Diffusion and aggregation} 

\abstract{
We describe a new phenomenon in models of coalescence and fragmentation, that of {\em gel-shatter cycles}. These are dynamical, unforced, stochastic cycles in which slow, approximately deterministic coalescence up to and beyond gelation is followed by abrupt random shattering. We describe their appearance in simulations of stochastic models with multiplicative kernels for coalescence and spontaneous fragmentation into monomers (`shattering'). The regime in which such cycles occur is characterized by a cyclicity order parameter, and we provide a simple scaling argument which describes both this regime and those which border it.
}

\begin{document}

\maketitle

\begin{figure}
    \onefigure[width=0.45\textwidth]{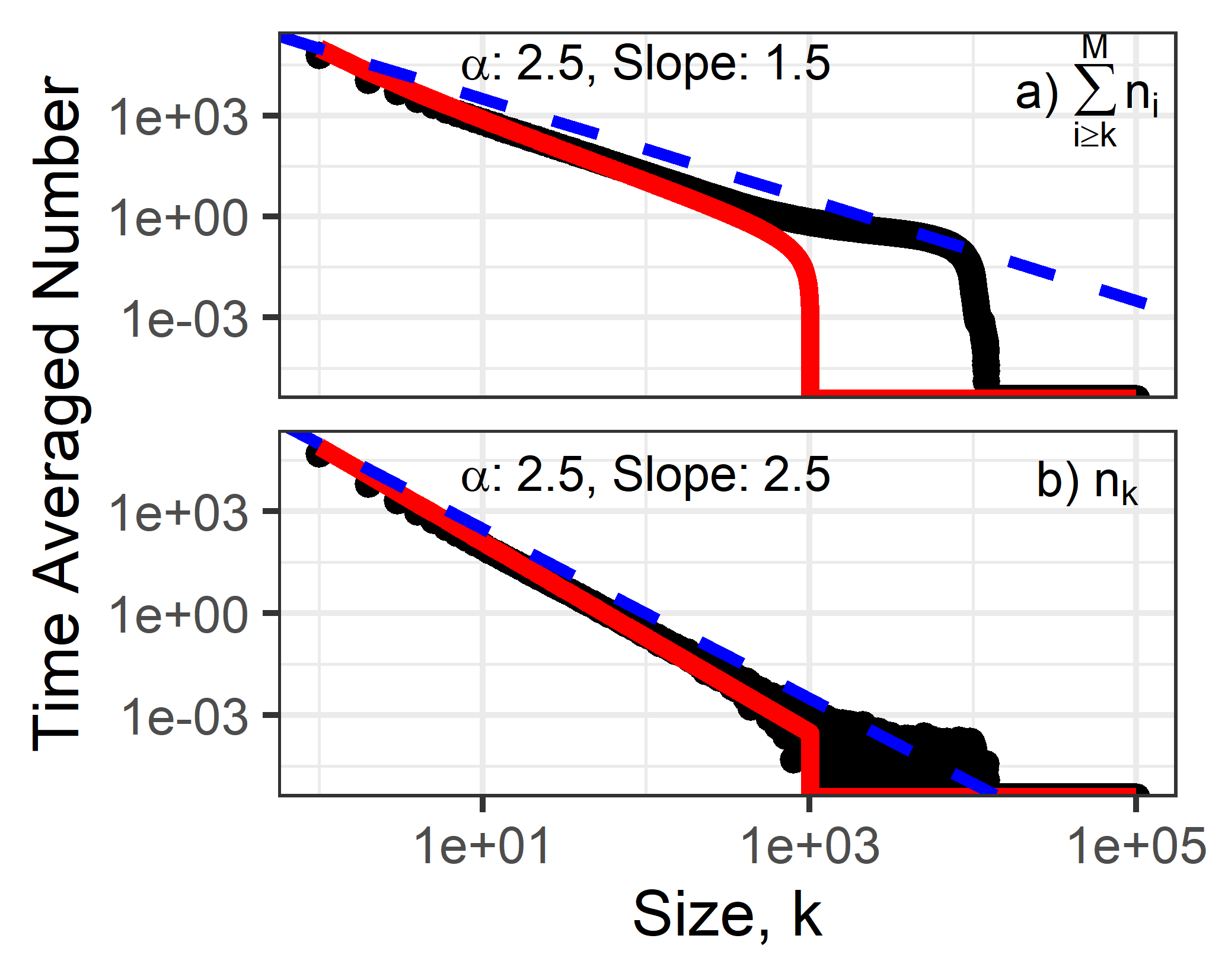}
    \caption{Cluster size distribution. A typical time-averaged distribution (black dots) emerging as a result of multiplicative coalescence and multiplicative shattering, with fragmentation restricted to clusters larger than $10^4$. The total mass is $M=10^5$. (a) Complementary cumulative distribution of clusters, {\em i.e.}\ the number of clusters larger than a certain size. (b) Mean cluster density. This follows a truncated $\alpha = \frac{5}{2}$ power-law (blue dashed line). We also show the distribution at a single time whose power-law exponent is closest to the median (red line).}
    \vspace{1cm}
    \label{fig:TimeAverage_Barrier}
\end{figure}

It is natural that coalescence of smaller units should result in larger accumulations that may eventually become unstable and fragment.
A wide range of models for such coalescence and fragmentation processes have been created, often with overlapping or even contradictory results, in fields as diverse as physical chemistry \cite{Spicer96_FlocCoFr}, probability \cite{Aldous99_CoalProbabilists}, fluid dynamics \cite{Pushkin02_Self}, social grouping \cite{Johnson05_Original}, bank mergers \cite{Pushkin04_BankMergers} and terrorist grouping networks \cite{Bohorquez09_Common, Clauset10_Wiegel, Johnson16_NewOnline}. 
Part of this popularity is undoubtedly due to the emergence of power-law behaviour from simple coalescence, \textit{e.g.\ }two asteroids collide \cite{Tanaka96_SelfSimilar}, one fish eats another \cite{Datta11_MarineCoFr} or two investors trade information \cite{DHulst00_EconoCoFr}, all of which yield power-law-distributed steady states, \textit{i.e.}\ $p(x) \propto x^{-\alpha}$ for $\alpha > 1$ (see Figure \ref{fig:TimeAverage_Barrier}). 
When combined, coalescence and fragmentation are expected to balance to give a non-trivial equilibrium \cite{Tanaka96_SelfSimilar,Ruszczycki09,Clauset10_Wiegel,Birnstiel11_Astrophysics,Connaughton18_CoShattering,Kyprianou18_Universality}.

\begin{widetext}

\begin{equation} \label{eq:GenCoFr}
    \frac{dn_k}{dt} = \frac{1}{2} \sum_{i = 1} ^ {k-1} K(i, (k - i)) n_i n_{k-i} - n_k\sum_{i = 1} ^ \infty  K(i, k) n_i  -  F(k) n_k + \sum_{i = k+1}^{\infty} F(i) b(i, k) \frac{i}{k} n_i \mbox{.}
\end{equation}

\end{widetext}

However, it has recently been discovered that the cluster-size distributions emerging from the tension between coalescence and fragmentation do not necessarily tend to a steady state. 
Perpetual oscillations were first observed in coalescence systems with constant influx of monomers and removal of the largest clusters if the coalescence rate between the large and small clusters increases sufficiently fast as a function of cluster sizes \cite{Ball12_CoalOscillations}. 
It was shown that here the transition from a steady state to an oscillatory state occurs via a Hopf bifurcation. 
Subsequently, running waves have been observed in deterministic mean-field numerical solutions and analytics for a special choice of coalescence and collision-controlled fragmentation rules in the absence of cluster influx/outflux \cite{Matveev17_OscillationsCoFr}. 
Most recently, it has been shown that deterministic temporal oscillations can arise in a class of coalescence-fragmentation (C-F) processes where clusters grow or shrink by addition or deletion of monomers \cite{Pego20_BeckDoerCycles}. 
All of these findings predicted deterministic oscillations and relied on mean-field descriptions.

Here we use simulations to show that stochastic cyclical dynamics -- in the form of steady accumulation followed by sudden shattering -- emerge in an important class of C-F processes. 
We consider spontaneous fragmentation, not induced by collision. Such C-F processes have been used to study polymerisation/depolymerisation reactions \cite{blatz1945note}, cloud droplets spectra \cite{pruppacher2010microstructure}, and assembly/disassembly of cell membrane microdomains \cite{richardson2007toward}, among other applications. We focus on a particular example of size-biased coalescence and fragmentation processes which is widely used in applications. We provide simple scaling arguments to explain the observed dependence of the recurrence time on the system size.

Mathematical descriptions of the problem can be grouped into two broad classes: mean field approaches and stochastic models (including stochastic simulations).
Coalescence processes, also known as coagulation or aggregation processes, were introduced to mathematics a century ago \cite{Smoluchowski16_Coagulation}, while the stochastic equivalent, the Marcus-Lushnikov process, arose in the 1960s \cite{Marcus68_Co,Lushnikov78_Co}. 

\begin{figure}[hb!]
    \includegraphics[width=0.45\textwidth]{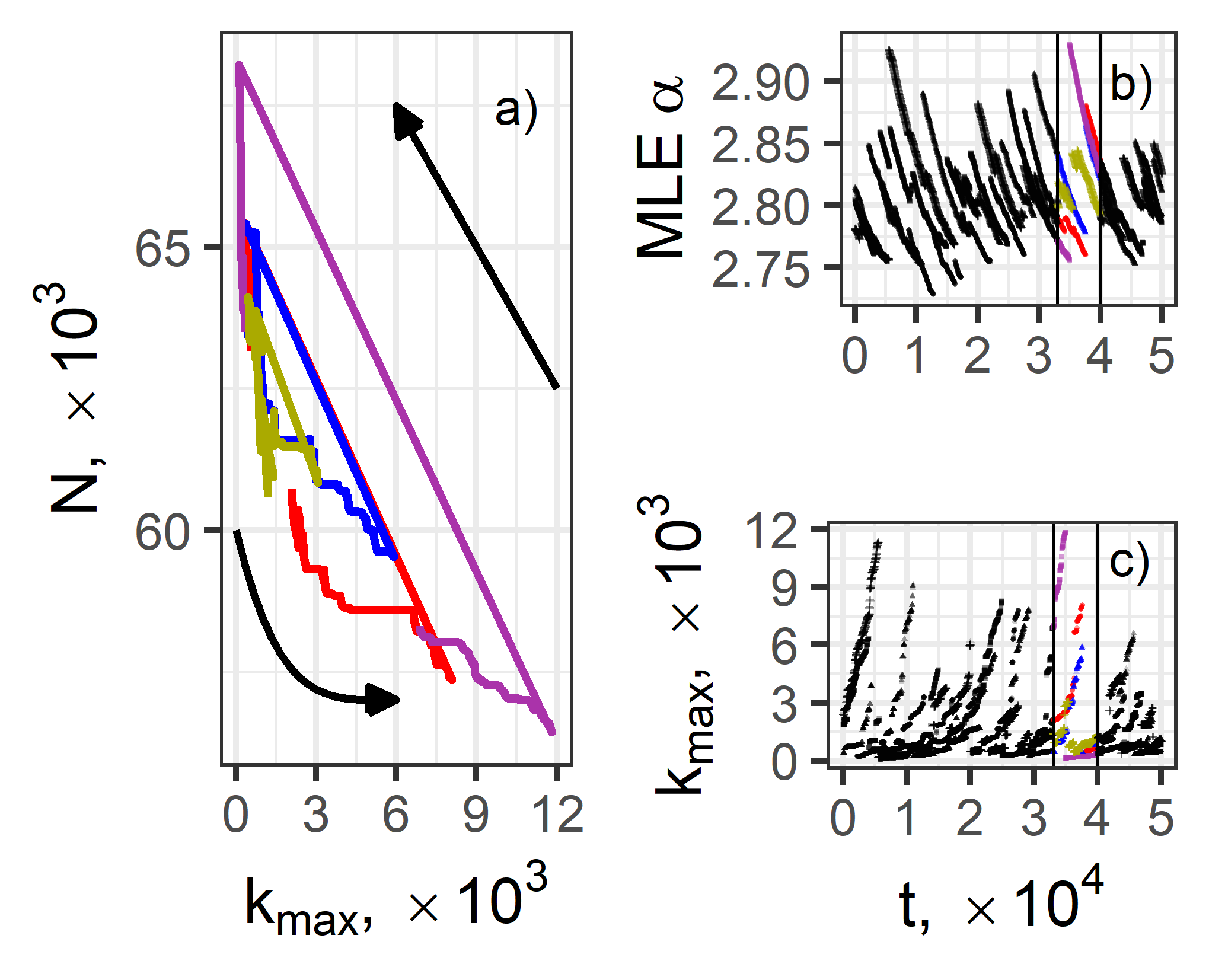}
    \caption{Stochastic coalescence-fragmentation cycles. 
    (a)~Number of clusters $N$ versus maximum cluster size $k_{max}$. 
    The counterclockwise trajectory shows periods of relatively gradual growth of the largest cluster followed by its abrupt shattering. 
    (b) The time-dependent power-law exponent $\alpha$ obtained by fitting the instantaneous cluster-size distribution to a (truncated) power-law using maximum likelihood estimation. The coloured region corresponds to the trajectories shown in (a). 
    (c) Maximum cluster size as a function of time.   
    Here $\hat K=0.99$, $\hat F=0.01$ and the total populations is $M=10^5$.
    An extended and animated version of this figure is included with the Supplemental Material.}
    \label{fig:Summary_Statistics_NoBarrier}
\end{figure}

First, for reference and convenience, we present an adaptation of Smoluchowski's equations for coalescence and fragmentation, modified from reference \cite{Spouge84_ExistenceCoFr}.
In equation (\ref{eq:GenCoFr}), $t$ is time and $n_k$ is the mean field density of clusters of (discrete) size $k$. 
Larger clusters are generated when two smaller clusters of sizes $i$ and $j$ randomly coalesce, in proportion to a kernel $K(i,j)$, to form a cluster of size $i + j$.
Clusters are reduced in size by a separate process of fragmentation, which occurs with the rate $F(i)$ for a cluster of size $i$. The function $b(i,j)$ describes what proportion of mass of the cluster of size $i$ goes to a cluster of size $j$; it is normalised so that $\sum_{j=1}^{i-1} b(i,j) = 1$. This condition guarantees that the mass is conserved in a fragmentation event. Since here fragmentation is spontaneous rather than collision-induced, the terms describing fragmentation are linear in the density of clusters.

The mean field approaches provide the basic language for cluster-size distribution modelling but ignore finite system size. However, in any computational implementation finite-size effects associated with the total system size $M$ will be important. In particular, the dynamics of pure coalescence, when the coalescence rate grows sufficiently fast as a function of the sizes of the coalescing clusters, leads to a giant cluster of a size comparable to $M$, a {\em gel} \cite{Aldous99_CoalProbabilists,Wattis06,Ziff82_GelationExponent, Lushnikov06_GelationReview}. \cite{Ziff82_GelationExponent}. In the thermodynamic limit, when $M \to \infty$ and the gel size is infinite, gelation goes beyond mean field theory (\ref{eq:GenCoFr}) and stochastic modelling is essential for capturing the ensuing dynamics. 

A particularly important coalescence kernel allowing gel formation naturally emerges in describing random network growth. When two nodes are randomly connected by an edge per unit time, two clusters of connected nodes coalesce at the rate $K = \hat K (i/M) (j/M)$, where $\hat K$ is constant, the system size $M$ is the total number of nodes, and $i/M$ and $j/M$ are the probabilities of picking clusters of sizes $i$ and $j$ from the population. This kernel is often referred to as the multiplicative coalescence kernel.

Similarly, the multiplicative fragmentation rate $F(i)$ is proportional to the cluster size $i$. The particular case in which the fragmented clusters disaggregate into monomers is referred to as {\em shattering} \cite{Birnstiel11_Astrophysics, Connaughton18_CoShattering} (We adopt the terminology used by reference \cite{Birnstiel11_Astrophysics}, but note that shattering has been used to describe alternative phenomena, \textit{e.g.\ }\cite{McGrady87_Shattering}.)
This form of fragmentation can be conveniently simulated numerically by picking a random node from the population of all nodes and shattering the cluster it belongs to, removing all edges amongst the formerly connected nodes. The corresponding fragmentation rule is given by
$F(i)=\hat F (i/M),\,  b(i,j)=\delta_{j 1}$, where $\hat F$ is constant and $i/M$ is the probability of picking a cluster of size $i$.

Simple kernels, especially those that are multiplicative, often result in solvable steady state systems, and the above forms are no different.
Under the assumption of a steady state, ${dn_k}/{dt} = 0$ for all $k$, set $\rho_k = k n_k / M$ so that for $k \geq 2$
\begin{equation*}
    \rho_k = \frac{\frac{1}{2} \hat K}{\hat F + \hat K \sum_{i = 1}^{\infty} \rho_i} \sum_{i = 1}^{k - 1} \rho_i \rho_{k - i}.
\end{equation*}
Provided that $\sum_{i = 1}^{\infty} \rho_i$ is well-behaved, which empirically is the case when $\hat F$ is sufficiently large compared to $\hat K$, we have a variation of the standard recurrence relation of the Catalan numbers \cite{Hilton91_Catalan,Catalan_OEIS}.
One standard method of solving such a system is to use generating functions and series expansion \cite{DHulst00_EconoCoFr,Clauset10_Wiegel}. 
Writing the result in terms of the Catalan numbers $C_n$ and writing the prefactor as $\gamma$, we have
\begin{equation}\label{eq:cofr:catsoln}
    \rho_k = C_{k - 1} \gamma^{k - 1} \rho_1^k.
\end{equation}
Stirling's law then yields the standard result of a truncated power-law with exponent $-3/2$.
The $\rho_1$ term can be evaluated by considering the total amount of mass in the system to be fixed.
Many other, related solutions exist for the steady-state distribution, dependent on the precise formulation of Equation \ref{eq:GenCoFr} \cite{Pushkin02_Self,Ruszczycki09,Birnstiel11_Astrophysics,Kyprianou18_Universality}.

The typical stochastic dynamics and its emerging cycles are seen in Figure \ref{fig:Summary_Statistics_NoBarrier}(a), which shows a sample of the trajectory in the $(k_{max},N)$ plane, where $N$ is the total number of clusters and $k_{max}$ is the size of the largest cluster. 
Each distinct counterclockwise cycle begins in the top-left with a predominance of monomers. 
A period of gradual growth of $k_{max}$ is accompanied by a decrease in $N$, taking the trajectory slowly down and to the right before accelerating rightward through gelation until growth ends in an abrupt random jump back to top-left with $\Delta N \approx - \Delta k_{max}$. 
This clearly evidences that 1) coalescence leads to the formation of very large clusters during the periods of gradual cluster size growth and 2) such periods end in shattering of the largest cluster. 
However, the cycling dynamics are not limited to the growth and shattering of the largest cluster but involves the whole cluster size distribution. Indeed, our simulations show that the cluster size distribution is very broad and can be fitted to (truncated) power-laws with the time-dependent exponent $\alpha(t)$ cycling in the range $2.7-2.9$. 
During the period of coalescence, the cluster size distribution is broadening and $\alpha$ decreases appreciably. 

To the best of our knowledge, such stochastic {\em gel-shatter cycles} have not been previously described. We hypothesize that they should be a generic feature of C-F systems with gel-forming coalescence and sufficiently strong fragmentation. 

\begin{figure}
    \centering
    \includegraphics[width=0.475\textwidth]{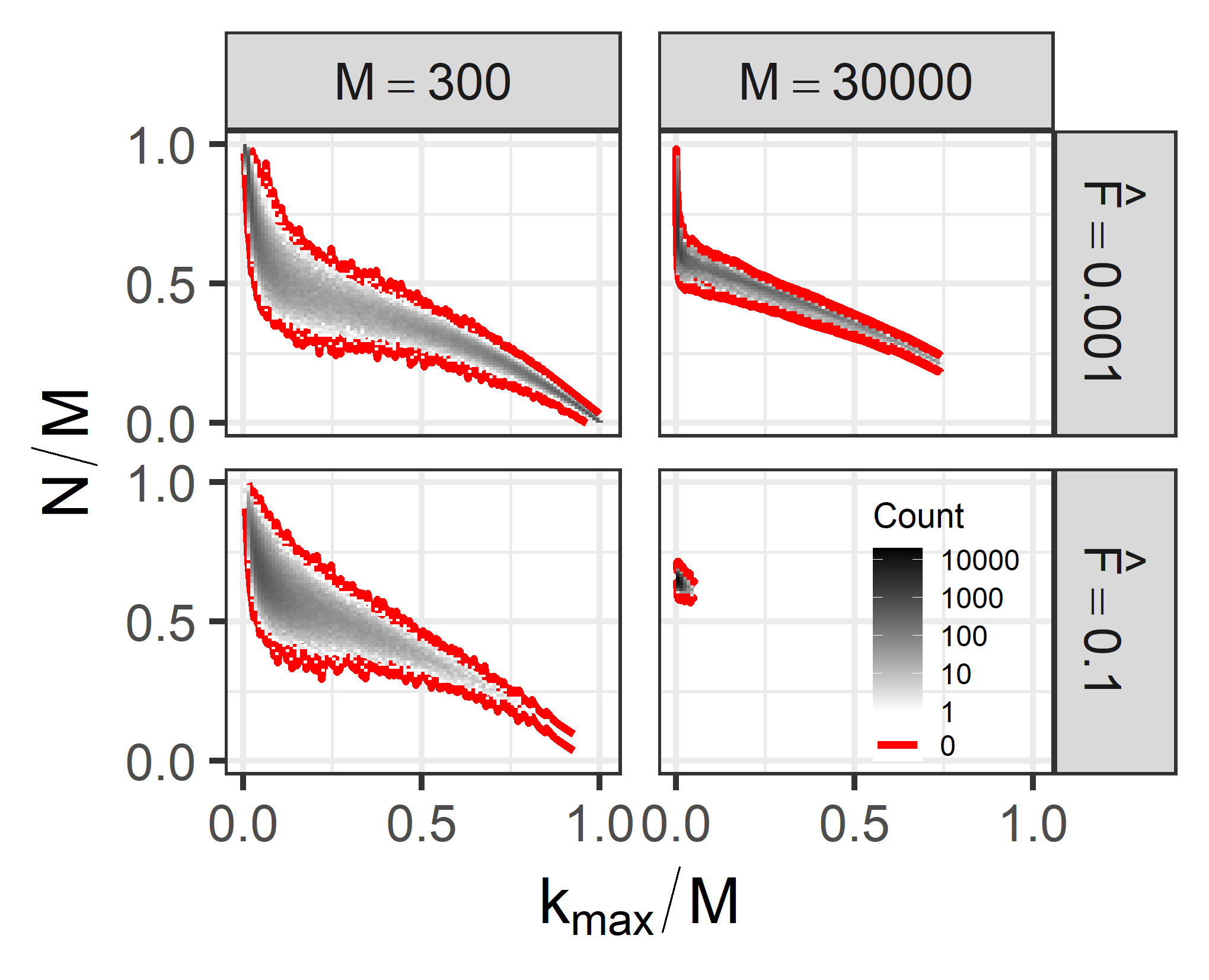}
    \caption{Cluster summary statistics. System sizes are $M = 3\times 10^2$ and $3\times 10^4$ (left to right), and fragmentation rates ${\hat F}= 10^{-3}$ and $10^{-1}$ (top to bottom, with coalescence rates ${\hat K}=1-{\hat F}$). The number of clusters $N$ and the maximum cluster size $k_{max}$ are normalised by $M$.  Darker regions represent states that emerge more often. The red border line denotes the boundary of the region visited by the system.
    A more detailed figure is present in the Supplemental Material.}

    \label{fig:HeatMaps_NoBarrier}
\end{figure} 

The gel-shatter cycles are strongly dependent on the number of monomers at complete disaggregation, which we term the system size $M$. 
Figure \ref{fig:HeatMaps_NoBarrier} presents a `heat map' of the times the simulation spends in different regions of the $(k_{max},N)$ plane (normalised by $M$ for comparison between different populations). 
For low fragmentation rates and small system sizes (top-left corner) the system visits a broad shoulder-like region. The `elbow', the change in gradient of its lower boundary, corresponds roughly to gelation.
As fragmentation rate and system size increase, the broad distribution collapses to a very small region within which the stochasticity is not visible as distinct cycles. 

\begin{figure*}
    \centering
    \includegraphics[width=0.95\textwidth]{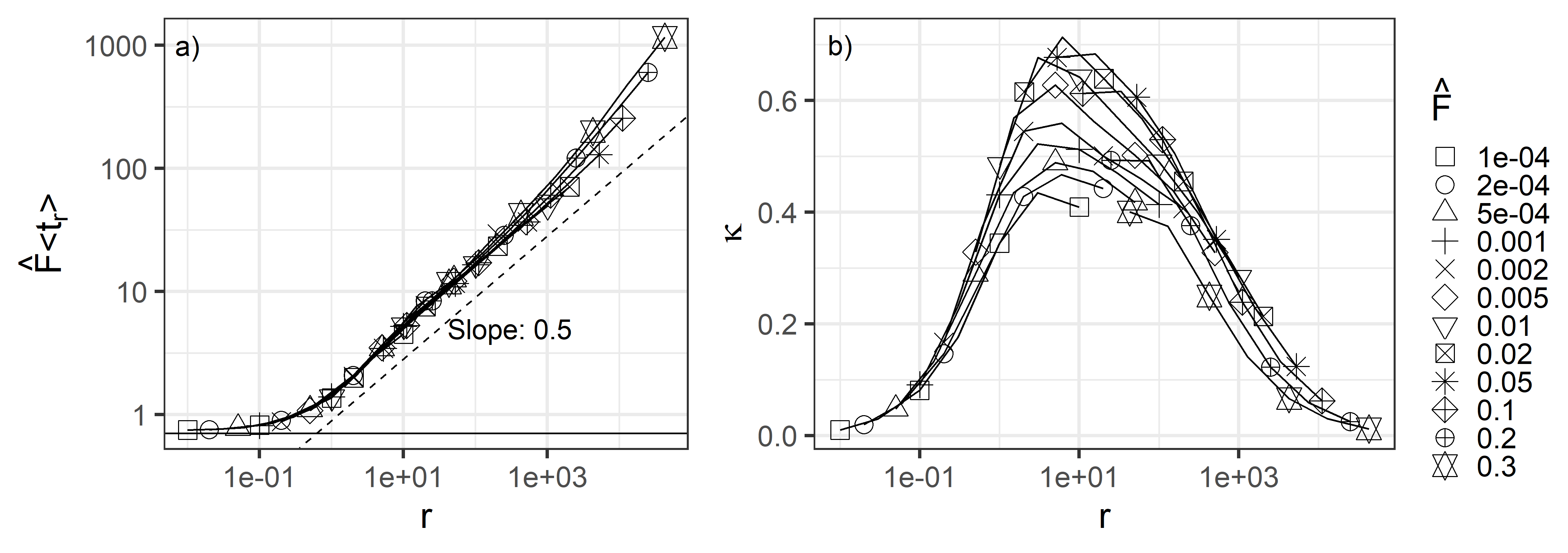}
    \caption{Observed data collapse in terms of $r$. Points are plotted for $M = 10^2, 10^3, 10^4$ and $10^5$, but data were additionally gathered at $M = 3 \times 10^2, 3 \times 10^3$ and $3 \times 10^4$.
    (a) When plotting the mean recurrence time $\meantr$ multiplied by the fragmentation rate $\hat F$, we see that the middle region of unforced gel-shatter cycles has strong data collapse following a nearly linear trend on log-log axes. 
    (b) When plotting the order parameter $\mathcal{K}$, we observe a plateau of $\mathcal{K}$ in the region of unforced gel-shatter cycles.}
    \label{fig:meanrunlength}
\end{figure*}

The mean recurrence time $\meantr$ is defined to be the number of computational steps between successive shatterings of the largest cluster, averaged across simulations with the same parameters. This can be thought of as the average duration of a cycle in the gel-shatter regime.
The dependence of the mean recurrence time $\meantr$ on system size $M$ for varying fragmentation rates $\hat F$ shows crossovers between three distinct regimes: 
    for small systems and small fragmentation rates, $\meantr \sim M^{0}$ 
    while modestly increasing $M$ or $\hat F$ yields ${\meantr \sim M^{1/2}}$,
    before crossing over into a regime with superlinear growth of $\meantr$ with $M$.
It is {\em a priori} possible that these regimes are due to finite size effects. In our simulations, finite-size effects are strong and appear even for systems of size $M = 10^5$.
Furthermore, there are multiple time scales present, and the importance of these could vary across rates and system sizes.
These curves show data collapse when $\hat F \meantr$ is plotted as a function of the dimensionless parameter
\begin{equation}
    r=\frac{\hat F M}{\hat K}.
\end{equation}
The physical meaning of $r$ is the ratio of the characteristic times of gelation $T_g \sim M/ \hat K$ and shattering $T_f \sim 1/\hat F$. 
Figure \ref{fig:meanrunlength}(a) demonstrates that 
\begin{equation}
    \langle t_r \rangle  = \hat F^{-1} g(r),
\end{equation}
where the scaling function $g(r) \sim 1$ for $r \ll 0.1$. 
For larger values of $r$, the scaling function crosses over to $g(r) \sim r^{1/2}$ and this behavior persists for almost four orders of magnitude of $r$. 
Finally, for $r \gtrsim 10^3$, the scaling breaks down. 

In order to explore further the nature of the stochastic gel-shatter cycles and distinguish cyclical dynamics from acyclical stochastic fluctuations, we introduce a cyclicity order parameter, $\mathcal{K}$, defined to be the number of computational steps that result in growth of the maximum cluster size $k_{max}$ minus those that result in reduction, divided by the total number of computational steps. 
Hence, $\mathcal{K}$ must lie in the range $-1$ to $1$.
A value near zero would indicate that the largest cluster experiences approximately as many coalescence events as it does fragmentation events, regardless of the magnitude of the events, and is characteristic of generic stochastic fluctuations which have no preferred sense.
A high value of $\mathcal{K}$, in contrast, indicates many steps of growth followed by abrupt collapse, thereby characterizing gel-shatter cycles.  
(Large negative $\mathcal{K}$ would indicate many fragmentation events followed by rare coalescence. This is possible only if fragmentation events are frequent and small, and we do not expect to observe it with our chosen $\hat K$ and shattering fragmentation kernel $F$.)
The order parameter $\mathcal{K}$ depends non-trivially on $r$, see Figure \ref{fig:meanrunlength}(b).
Here $\mathcal{K}(r)$ forms hump-shaped curves with the maximum around $r=10$. 
The order parameter is appreciably large, $\mathcal{K}>0.1$, for $0.1<r<10^3$, which constitutes the middle regime with ${\langle t_r \rangle \sim M^{1/2}}$.

We can now identify the different scalings as three physical regimes:
\begin{enumerate}[(i)]
    \item {\bf Weak fragmentation or forced cycles.} Here ${T_g \ll T_f}$ and coalescence quickly results in a single gel cluster. The gel is shattered at the rate $\hat F$ and hence the emerging cycles have the mean recurrence time  $\langle t_r \rangle \sim \hat F^{-1}$. This regime is physically unsurprising. 
    
    \item  {\bf Unforced gel-shatter cycles.} Here $T_g \sim T_f$ and the cycles arise from a non-trivial interplay between the dynamics of gelation and shattering. The cyclicity  $\mathcal{K}$ reaches a maximum in this regime. The presence and extent of this regime in gelling C-F systems is our main finding.
    
    \item {\bf Fragmentation-dominance.}  Here $T_g \gg T_f$ and fragmentation is so strong as to preclude formation of large clusters. This annihilates the gel-shattering cycles, leaving only stochastic variation in a small region.
\end{enumerate}
    
We now provide a simple scaling argument to explain the behaviour of the mean recurrence times and gain further insight into the gel-shatter cycles.
    
Assume the dynamics are dominated by continuous growth and shattering of the largest cluster. 
Writing the size of the largest cluster $k_{\max}$ at time $t$ as $m(t)$, we track the largest cluster and the probability, dependent on $m(t)$, that it shatters at a given time. 
The probability $P(t)$ that the largest cluster is not shattered by time $t$ is governed by  
\begin{equation}
    \frac{dP(t)}{dt} = - \hat F M^{-1} m(t) P.
\end{equation}
Hence the probability density for shattering at time $t$ is 
\begin{equation}
p(t)= - \frac{dP(t)}{dt} = \hat FM^{-1} m(t) e^{ -\hat F M^{-1} \int_0^t m(\tau) d\tau}.
\end{equation}
In the forced-cycles regime (i), coalescence quickly leads to a single large cluster of size $M$, hence $m(t)=M$ and
\begin{equation*}
    p(t) = \hat F  \exp{ \left( - \hat F  t \right) }.
\end{equation*}
This is the exponential distribution with mean ${\meantr = \hat F^{-1}}$; thus  $g(r)\sim 1, \, r \ll 1$.
   
In the unforced gel-shattering cycles regime (ii), the simplest form for $m(t)$ consistent with dimensional arguments is linearity. We assume that $m(t) = c \hat K t$, where $c>0$. Then
$$
    p(t) = c \hat F \hat K M^{-1} t  e^{  - c \hat F \hat K M^{-1}  t^2/2  }.
$$
This is the Rayleigh distribution with scale parameter $\left(\sqrt{c \hat F \hat K M^{-1}}\right)^{-1}$ and mean
\begin{equation}
    \meantr = \sqrt{\pi M/(2 c \hat F \hat K)} .
\end{equation}
Hence, $g(r) \sim r^{1/2}$ in regime (ii). 
Also, it is straightforward to show that the largest realisable cluster size scales as
$M r^{-1/2}$.
This scaling explains the shortening of the `elbow' on Figure \ref{fig:HeatMaps_NoBarrier} for larger $M$ and $\hat F$.
    
This work reports a distinctive new phenomenon, that of stochastic {\em gel-shatter cycles}. 
Gel-shatter cycles are observed explicitly in simulations of multiplicative (size-biased) coalescence and spontaneous-shattering fragmentation kernels when the time scales for these are in balance.
We expect gel-shatter cycles to be a general phenomenon beyond that of multiplicative kernels.
The prerequisite qualities appear to be three-fold: gelation, which promotes the growth of a single cluster that dominates the system; a strong form of fragmentation, which can reset the system to a pre-gel state; and balanced time scales, so that the gel is neither instantly removed nor aggregates the entire system.
Gelation is a broad and well-studied phenomenon, occurring with many coalescence kernels.
Shattering is perhaps the strongest assumption we make, but we anticipate that this can be relaxed to a sufficiently strong form of mixed fragmentation. 
Finally, in any system studied for its tension between coalescence and fragmentation (as opposed to just one of these alone), the time scales of these are likely to be comparable.
In such circumstances practitioners should be wary of presuming the existence of a steady state, and conscious that gel-shatter cycles may be present. 

On this basis we anticipate that gel-shatter cycles are a ubiquitous phenomenon.
This salient novel feature emerges from a range of simulation studies originally intended to explore the robustness of the Smoluchowski coalescence-fragmentation model and its variants (see Supplemental Material). 
Whilst our conclusion is that models of this type are robust to perturbations of various kinds -- confirming the common view \cite{Bohorquez09_Common,Johnson16_NewOnline,Ruszczycki09,Clauset10_Wiegel} -- we note the ubiquitous presence of unforced stochastic gel-shatter cycles in finite-size systems. 
This is relevant to the application of these models to physical and social phenomena, where the gel-shatter cyclicity may be a natural dynamical feature, as well as affecting the robustness of the computationally-fitted exponents.
The presence of a non-trivial dynamic underlying a presumed steady state needs to be considered when models of this type are used in applications, and is an intriguing avenue for future theoretical research.

\acknowledgments
We would like to thank Stephen Connor and Gustav Delius for interesting discussions. We are also grateful to two anonymous referees for their constructive comments.

\bibliographystyle{eplbib}
\bibliography{references}
\end{document}